\documentclass{amsart}
\usepackage{amssymb}
\usepackage{a4}
\begin{document}
\newtheorem{theorem}{Theorem}
\def\qedbox{\hbox{$\rlap{$\sqcap$}\sqcup$}}
\def\a{\alpha}
\def\b{\beta}
\def\c{\gamma}
\def\la{\langle}
\def\ra{\rangle}

\def\w{\mathcal W}
\def\J{\mathcal J}
\def\C{\mathbb{C}}

%\makeatletter
% \renewcommand{\theequation}{%
% \thesection.\alph{equation}}
% \@addtoreset{equation}{section}
% \makeatother
\title[Conformally Osserman manifolds and self-duality]
{Conformally Osserman manifolds and self-duality in Riemannian
geometry}
\author{Novica  Bla{\v z}i{\'c} and  Peter Gilkey}
\begin{address}{NB: Faculty of Mathematics, University of Beograd, Studentski Trg. 16, P.P. 550,
11000 Beograd, Srbija i Crna Gora. Email: {\it blazicn@matf.bg.ac.yu}}\end{address}
\begin{address}{PG: Mathematics Department, University of Oregon,
Eugene Or 97403 USA.\newline Email: {\it gilkey@darkwing.uoregon.edu}}
\end{address}

\begin{abstract}
We study the spectral geometry of the conformal Jacobi operator on a $4$-dimensional
Riemannian manifold $(M,g)$. We show that $(M,g)$ is conformally Osserman if and only
if $(M,g)$ is self-dual or anti self-dual. Equivalently, this means that the
curvature tensor of $(M,g)$ is given by a quaternionic structure, at least pointwise.
\end{abstract}
%\version}
\maketitle

 Let $\mathcal{R}_g$ and $\mathcal{W}_g$ be
the curvature operator and the Weyl curvature operator associated to a Riemannian
manifold $(M,g)$. The {\it Jacobi operator} $\mathcal{J}_{\mathcal{R}_g}$ and the {\it
conformal Jacobi operator}
$\mathcal{J}_{\mathcal{W}_g}$ are defined by:
$$
\mathcal{J}_{\mathcal{R}_g}(x):y\rightarrow\mathcal{R}_g(y,x)x\quad\text{and}\quad
\mathcal{J}_{\mathcal{W}_g}(x):y\rightarrow
\mathcal{W}_g(y,x)x\,.
$$
In contrast to the Jacobi operator, the {\it conformal Jacobi operator} is
conformally invariant \cite{BG}; if $h=e^\alpha g$ is a conformally
equivalent Riemannian metric, then
$$\mathcal{J}_{W_h}=\mathcal{J}_{W_g}\,.$$

One says that $(M,g)$ is {\it Osserman} (resp. {\it conformally Osserman}) if the
eigenvalues of $\mathcal{J}_{R_g}$ (resp. $\mathcal{J}_{W_g}$) are constant on the
bundle of unit tangent directions. In a series of papers started by Chi \cite{Ch88}
and continued by Nikolayevsky
\cite{Nik,Nik2,Nik3} it was shown that Osserman manifolds of dimension $n\neq 16$
are two-point homogeneous spaces; for $m\ne16$, this gives an affirmative answer to
a question raised by Osserman (see, for example, \cite{Oss}).

Previous work \cite{BG} has shown that conformally Osserman manifolds are
conformally flat if $m\equiv1$ mod $2$ and are either conformally flat or conformally
equivalent to a complex space form (i.e. to complex
projective space with the Fubini-Study metric or to the non-compact dual) if
$m\equiv2$ mod
$4$. In this present paper, we study the $4$-dimensional setting motivated by work
of Sekigawa and Vanhecke (\cite{SV}). Our main result is the following.

\begin{theorem}\label{thm-1} Let $(M,g)$ be a $4$-dimensional oriented Riemannian
manifold. The following conditions are equivalent:
\begin{enumerate}
\item $(M,g)$ is conformally Osserman.
\item $(M,g)$ is self-dual or anti-self dual.
\end{enumerate}
\end{theorem}

It is useful to work in a purely algebraic setting. Let $V$ be a finite
dimensional vector space equipped with a positive definite inner product
$(\cdot,\cdot)$. We say that $R\in\otimes^4V^*$ is
an {\it algebraic} curvature tensor if $R$ has the usual symmetries:
\begin{eqnarray*}
&&R(x,y,z,w)=-R(y,x,z,w)=R(z,w,x,y),\quad\text{and}\\
&&R(x,y,z,w)+R(y,z,x,w)+R(z,x,y,w)=0\,.
\end{eqnarray*}
The associated curvature operator $\mathcal{R}(x,y)$ is then characterized by the
identity
$(\mathcal{R}(x,y)z,w)=R(x,y,z,w)$. Let $\tau$ be the scalar curvature and let
$\rho$ be the Ricci operator. Set
\begin{eqnarray*}
&&\mathcal{R}_g^0(x,y)z:=g(y,z)x-g(x,z)y,\quad\text{and}\\
&&\mathcal{L}(x,y)z:=g(\rho y,z)x-g(\rho x,z)y+g(y,z)\rho x-g(x,z)\rho y\,.
\end{eqnarray*} 
The associated Weyl operator given by
\begin{eqnarray*}
&&\mathcal{W}(x,y):=\mathcal{R}(x,y)
   +\frac1{(m-1)(m-2)}\tau\mathcal{R}_g^0(x,y)+\frac1{m-2}\mathcal{L}(x,y)
\end{eqnarray*}

We say that $R$ is {\it Osserman} (resp. {\it conformally Osserman}) if the
eigenvalues of
$\mathcal{J}_R$ (resp. $\mathcal{J}_W$) are constant on the unit sphere in $V$.

Let $\Phi$ be a skew-symmetric endomorphism of $V$ with
$\Phi^2=-1$. Define
$$R_\Phi(x,y)z:=(\Phi y,z)\Phi x-(\Phi x,z)\Phi y-2g(\Phi x,y)\Phi z\,.$$
Then $R_\Phi$ is an {\it algebraic curvature tensor}. We say that
$\{\Phi_1,\Phi_2,\Phi_3\}$ is a {\it unitary quaternion} structure on $V$ if the
$\Phi_i$ are self-adjoint and if the usual structure equations are satisfied:
$$\Phi_i\Phi_j+\Phi_j\Phi_i=-2\delta_{ij}\operatorname{id}\,.$$

Theorem \ref{thm-1} is
a consequence of the following purely algebraic fact:

\begin{theorem}\label{thm-2} Let $R$ be an algebraic curvature tensor on a
$4$-dimensional vector space $V$. The following assertions are equivalent:
\begin{enumerate}
\item $R$ is conformally Osserman.
\item $R$ is self-dual or anti-self dual.
\item There exists a unitary quaternion structure on $V$ so that\newline
$W=\lambda_1R_{\Phi_1}+\lambda_2R_{\Phi_2}+\lambda_3R_{\Phi_3}$
where $\lambda_1+\lambda_2+\lambda_3=0$.
\end{enumerate}
\end{theorem}

We begin by showing that the first assertion implies the second assertion in
Theorem
\ref{thm-2}.  Let
$R$ be a conformally Osserman algebraic curvature tensor on $\mathbb{R}^4$. Zero is
always an eigenvalue of $\mathcal{J}_W$ since $\mathcal{J}_W(x)x=0$. Let $e_1$ be a
unit vector. Since $\mathcal{J}_W(\cdot)$ is symmetric, it has an orthonormal
basis of real eigenvectors. Thus we may extend
$e_1$ to an orthonormal basis
$\{e_1,e_2,e_3,e_4\}$ so that
$$\mathcal{J}_W(e_1)e_2=ae_2,\quad
  \mathcal{J}_W(e_1)e_3=be_3,\quad
  \mathcal{J}_W(e_1)e_4=ce_4\,.$$
Since $\operatorname{Tr}(\mathcal{J}_W)=0$, we have $a+b+c=0$.

The argument given by Chi \cite{Ch88} in his analysis of the $4$-dimensional
setting was in part purely algebraic. This algebraic argument extends without
change to this setting to show that, after possibly replacing $\{e_2,e_3,e_4\}$ by
$\{-e_2,-e_3,-e_4\}$ that the non-vanishing components of the Weyl curvature on
this basis are given by:
\begin{equation}\label{eqn-1}\begin{array}{l}
W_{1221}=W_{3443}= - W_{1234}=a,\\
W_{1331}=W_{2442}= - W_{1342}=b,\vphantom{\vrule height 11pt}\\
W_{1441}=W_{2332}= - W_{1423}=c\,.\vphantom{\vrule height 11pt}
\end{array}\end{equation}
Let $e^{ij}:=e^i\wedge e^j$ where $\{e^i\}$ is the dual basis for $V^*$. We consider
the following bases for
$\Lambda_2^\pm(V)$:
$$
f_1^\pm=e^{12}\pm e^{34}, \ f_2^\pm=e^{13}\mp e^{24},\ f_3^\pm=e^{14}\pm
e^{23}\,.$$ Since $\w(e^{pq})=\frac{1}{2}W_{pqij}e^{ij}$,
$$\begin{array}{lll}
\mathcal{W}(f_1^-)=0,&\mathcal{W}(f_2^-)=0,&\mathcal{W}(f_3^-)=0,\\
\mathcal{W}(f_1^+)=-2af_1^+,&\mathcal{W}(f_2^+)=-2bf_2^+,&\mathcal{W}(f_3^+)=-2cf_3^+\,.
\end{array}$$
Thus conformally Osserman algebraic curvature tensors are self-dual.

Next we show that Assertion (2) implies Assertion (3) in Theorem
\ref{thm-2}. Suppose that $R$ is a self-dual algebraic curvature tensor
on $\mathbb{R}^4$. Let $e_1$ be a unit vector. Choose an orthonormal basis
$\{e_1,e_2,e_3,e_4\}$ for $\mathbb{R}^4$ so that
$$\mathcal{J}_W(e_1)e_2=ae_2,\quad\mathcal{J}_W(e_1)e_3=be_3,
\quad\mathcal{J}W(e_1)e_4=ce_4\,.$$
We then have
\begin{equation}\label{eqn-2}\begin{array}{lll}
W_{1221}=a,&W_{1231}=0,&W_{1241}=0,\\
W_{1321}=0,&W_{1331}=b,&W_{1341}=0,\\
W_{1421}=0,&W_{1431}=0,&W_{1441}=c\,.
\end{array}\end{equation}
By replacing $\{e_2,e_3,e_4\}$ by $\{-e_2,-e_3,-e_4\}$ if necessary, we can assume
that $\{e_1,e_2,e_3,e_4\}$ is an oriented orthonormal basis. We have
\begin{eqnarray*}
\mathcal{W}(e^{12}-e^{34})&=&(W_{1212}-W_{3412})e^{12}+(W_{1234}-W_{3434})e^{34},\\
&+&(W_{1213}-W_{3413})e^{13}+(W_{1214}-W_{3414})e^{14}\\
&+&(W_{1223}-W_{3423})e^{23}+(W_{1224}-W_{3424})e^{24}\,.
\end{eqnarray*}
Since $W$ is self-dual, $\mathcal{W}(e^{12}-e^{34})=0$.
Equation (\ref{eqn-2}) then implies
$$\begin{array}{lll}
W_{3412}=W_{1212}=-a,&
W_{3434}=W_{1234}=-a,&
W_{3413}=W_{1213}=0,\\
W_{3414}=W_{1214}=0,&
W_{3423}=W_{1223}=0,&
W_{3424}=W_{1224}=0\,.
\end{array}$$
We argue
similarly using
$e^{13}+e^{24}$ and
$e^{14}-e^{23}$ to see that the formulas of Equation (\ref{eqn-1}) hold. We define
a unitary quaternion structure by defining $\Phi_3:=\Phi_1\Phi_2$ where
$$\begin{array}{llll}
\Phi_1:e_1\rightarrow e_2,&\Phi_1:e_2\rightarrow -e_1,&
\Phi_1:e_3\rightarrow e_4,&\Phi_1:e_4\rightarrow-e_3,\\
\Phi_2:e_1\rightarrow e_3,& \Phi_2:e_3\rightarrow -e_1,&
\Phi_2:e_4\rightarrow -e_2,&\Phi_2:e_2\rightarrow e_4\,.
\end{array}$$
It is then immediate that the formulas of Equation (\ref{eqn-1}) hold for $\tilde
W:=aW_{\Phi_1}+bW_{\Phi_2}+cW_{\Phi_3}$
and thus $W=\tilde W$. This shows the second assertion implies the third assertion.

Finally, if $W$ is given by a unitary quaternion structure, then the discussion of
\cite{Gi02} shows that $W$ is Osserman. This completes the proof of Theorem
\ref{thm-2} and thereby of Theorem \ref{thm-1} as well. \hfill\qedbox.

\section*{Acknowledgments} 
Research of N. Bla{\v z}i{\'c} is partially supported by the  DAAD  (Germany) and MNTS
 Project 1854 (Serbia).
Research of P. Gilkey partially supported by the
Max Planck Institute in the Mathematical Sciences (Leipzig).
The first author wish to express his thank to the Technical University of Berlin where part of  research reported here was conducted.
Both authors wish to express their appreciation to the referee for helpful
suggestions concerning the matter at hand.

\end{document}